\documentclass{article}

\usepackage{arxiv}

\usepackage[utf8]{inputenc} 
\usepackage[T1]{fontenc}    
\usepackage{hyperref}       
\usepackage{url}            
\usepackage{booktabs}       
\usepackage{amsfonts}       
\usepackage{nicefrac}       
\usepackage{microtype}      
\usepackage{doi}
\usepackage{style}
  
\title{The Neural Approximated Virtual Element Method for Elasticity Problems}

\author{ \href{https://orcid.org/0000-0001-8642-4258}{\includegraphics[scale=0.06]{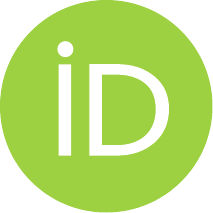}\hspace{1mm}Stefano~Berrone} \\
	Dipartimento di Scienze Matematiche\\
	``G. L. Lagrange''\\
	Politecnico di Torino, TO, 10129 \\
	\texttt{stefano.berrone@polito.it} \\
	\And
	\href{https://orcid.org/0000-0002-2837-2792}{\includegraphics[scale=0.06]{orcid.pdf}\hspace{1mm}Moreno~Pintore} \\
	Laboratoire Jacques-Louis Lions, \\
	Sorbonne Universit\'e, \\
        MEGAVOLT Team, Inria, \\
	4 place Jussieu, 75005 Paris, France \\
	\texttt{moreno.pintore@sorbonne-universite.fr} \\
	\And
	\href{https://orcid.org/0000-0002-8540-3639}{\includegraphics[scale=0.06]{orcid.pdf}\hspace{1mm}Gioana~Teora} \\
	Dipartimento di Scienze Matematiche\\
	``G. L. Lagrange''\\
	Politecnico di Torino, TO, 10129 \\
	\texttt{gioana.teora@polito.it} \\
}

\renewcommand{\shorttitle}{The Neural Approximated Virtual Element Method for Elasticity Problems}

\hypersetup{
pdftitle={The Neural Approximated Virtual Element Method for Elasticity Problems},
pdfsubject={q-bio.NC, q-bio.QM},
pdfauthor={Stefano~Berrone, Moreno~Pintore, Gioana~Teora},
pdfkeywords={NAVEM, Elasticity, Virtual Element Method, Neural Network, Basis Functions, Polygonal Meshes},
}
\begin{document}

\maketitle


\begin{abstract}
We present the Neural Approximated Virtual Element Method to numerically solve elasticity problems. This hybrid technique combines classical concepts from the Finite Element Method and the Virtual Element Method with recent advances in deep neural networks. Specifically, it is a polygonal method in which the virtual basis functions are element-wise approximated by a neural network, eliminating the need for stabilization or projection operators typical of the standard virtual element method. We present the discrete formulation of the problem and provide numerical tests on both linear and non-linear elasticity problems, demonstrating the advantages of having a simple discretization, particularly in handling non-linearities.
\end{abstract}


\section{Introduction}
Discretization methods relying on meshes comprising general polygons are receiving a great deal of attention from the scientific community. Indeed, polytopal methods can be very useful for a large variety of reasons, including the automatic handle of hanging nodes, the use of moving meshes, the ease of meshing complex geometries, and adaptivity \cite{Fassino2023, GrappeinTeora2025, Vicini2024}. Moreover, these methods have been successfully applied to several fields of computational mechanics, revealing several advantages over methods employing classical triangular/tetrahedral and quadrilateral elements. For instance, polygonal methods have shown great advantages in handling crack propagation and branching \cite{Leon2014}, contact problems \cite{Biabanaki2014}, and topology optimization \cite{Gain2015}.

In this framework, an important polygonal method is the Virtual Element Method (VEM), introduced in \cite{LBe13, LBe16} for elliptic problems and then extended to elastic problems in \cite{DaVeigaBrezzi2013, DaVeiga2015}. VEM can be seen as a generalization of the standard Finite Element Method (FEM) \cite{brennerscott, ciarlet2002finite} on more general meshes. One of the main features of the VEM, is that the local basis functions are not known in closed form. Indeed, even if they can be defined as the solution of local differential problems, they are never explicitly computed or approximated in the virtual element framework. It is therefore necessary to introduce suitable projection and stabilization operators to derive computable discrete bilinear forms associated with the Partial Differential Equation (PDE) of interest. In particular, the usage of a stabilization term is an undesirable aspect of the VEM method, as it lacks a theoretical foundation, remains problem-dependent, and an improper choice can significantly affect the method accuracy \cite{Teora2023}. Moreover, when dealing with non-linearities, the presence of the stabilization term makes the construction of the method more involved and increases the computational effort needed to solve the problem. Finally, the need for a projector to access the point-wise evaluation of virtual functions makes the computation of post-process quantities more difficult \cite{Chi2016}.

Because of these reasons, in the last few years, several numerical methods have emerged to address these issues while preserving the main advantages of the method. A non-exhaustive list of such methods is the following. In \cite{ChenSukumar2023, ChenSukumar2023b}, a stabilization-free formulation of the method, derived by enlarging the local projection polynomial order and enhancing the function spaces according to the polygon geometry, is presented. In \cite{Credali2024, TrezziZerbinati2024}, two approaches are proposed to locally approximate the local virtual basis functions and to get rid of any stabilization or projection operators: one based on the reduced basis method \cite{hesthaven2016certified} and the other on the \textit{Laplace solver} introduced in \cite{Gopal2019}. We observe that, in these last two methods, the basis functions are explicitly computed on each element, meaning that the ``virtual'' attribute only refers to the underlying discrete space.

Within the same context, the Neural Approximated Virtual Element Method (NAVEM) is introduced in \cite{PintoreTeora2024, PintoreTeora2025} for general elliptic problems. The main idea of this novel method is to use neural networks \cite{lecun2015deep} to accurately predict on the fly the local virtual basis functions, which are then used like in a standard FEM framework to compute a discrete solution. For this reason, the NAVEM method falls within the general recent framework of Scientific Machine Learning (SciML) \cite{Cuomo2022}, by combining the advantages of the non-linear approximation capabilities of neural networks with the ones of more classical numerical methods. This framework is exponentially increasing in popularity because it represents the intersection between two active communities that rarely interacted before. We highlight that several libraries are available to easily implement and adapt complex neural networks. The most used are Tensorflow \cite{tensorflow2015-whitepaper}, PyTorch \cite{pytorch}, and JAX \cite{jax2018github}.

In this paper, we present the NAVEM for linear and non-linear elastic problems. Moreover, we further highlight some new key aspects of the neural network that allow us to improve the NAVEM accuracy. It is very important to notice that the neural networks used in the NAVEM framework are trained using only the geometric information of the elements. Therefore, after a potentially expensive training phase, the same neural networks can be used multiple times to solve problems characterized by different material coefficients, constitutive laws, body forces, initial configurations but also, and most importantly, on various meshes. In this sense, it can be regarded as a method based on the so-called \textit{offline-online} splitting.

The manuscript is organized as follows. In Section \ref{sec:model_problem}, we present the model problem and its VEM discretization, describing its main properties, which are essential to devise a proper neural network training strategy, and highlighting its main issues. In Section \ref{sec:navem}, we introduce the NAVEM for elasticity problems. In particular, in Section \ref{sec:local_space}, we construct the function space from where the neural networks select the local basis functions in the NAVEM discretization, in Sections \ref{sec:nn_architecture} and \ref{sec:nn_training}, we describe the neural networks architecture and training, respectively, and in Section \ref{sec:navem_scheme} we discuss the general structure of the NAVEM as a hybrid method mixing ideas from the machine learning and scientific computing communities. Later, in Section \ref{sec:numerical_results}, we prove through numerical experiments the effectiveness of the method for elasticity problems, and the advantages of its simple formulation. The conclusion of the present work and the future perspectives are discussed in Section \ref{sec:conclusion}.

\section{The Model Problem and The Virtual Element Method}\label{sec:model_problem}

Let us consider a homogeneous isotropic elastic body $\Omega \subset \R^2$ with boundary $\Gamma$. We assume that this body is clamped on $\Gamma_D \subseteq \Gamma$, $\vert \Gamma_D \vert \neq 0$, and subjected to a body load $\ff \in \vleb{2}{\Omega}{2}$.

In the case of static problems, the governing equation of the solid in the initial configuration is
\begin{equation}
    \begin{cases}
        -\nabla \cdot \ssigma(\xx, \nabla \uu) = \ff(\xx) & \text{in } \Omega,\\
        \uu = \bm{0} & \text{on } \Gamma_D,\\
        \ssigma(\xx, \nabla \uu) \nn = \bm{0} & \text{on } \Gamma_N := \partial \Omega \setminus \Gamma_D,
    \end{cases}
    \label{eq:model_problem}
\end{equation}
where $\uu : \Omega \to \R^2$ represents the body displacement, whereas the second-order tensor $\ssigma = \ssigma(\xx, \nabla \uu)$ represents a constitutive law for the material at each point $\xx \in \Omega$.

Let us consider the space
\begin{equation*}
    \VM = \vsob[0,\Gamma_D]{1}{\Omega}{2} := \{ \vv \in \vsob{1}{\Omega}{2}:\ \vv = \bm{0} \text{ on } \Gamma_D\}
\end{equation*}
and the bilinear form
\begin{equation}
    \dbilin{\uu}{\vv} = \int_{\Omega} \ssigma(\xx, \nabla \uu (\xx)) : \nabla \vv,
    \label{eq:continuous_a}
\end{equation}
where, given two second-order tensors $\ssigma$ and $\Tmatrix$, we define $\ssigma : \Tmatrix = \sum_{i,j = 1}^n \ssigma_{ij} \Tmatrix_{ij}$.
The variational formulation of Problem \eqref{eq:model_problem} reads as: \textit{Find $\uu \in \VM$ such that}
\begin{equation}
    \dbilin{\uu}{\vv} = \int_{\Omega} \ff \cdot \vv \quad \forall \vv \in \VM.
    \label{eq:var_problem}
\end{equation}
We remark that here we are not considering internal constraints, such as incompressibility, which usually require suitable treatments. 

\subsection{The Virtual Element Space}

Let us consider a tessellation $\Th$ of the domain $\Omega$ made up of polygonal elements $E$. We denote by $h_E$, $\Nv[E]$, $\xx_E$ and $\Eh[E]$ the diameter, the number of vertices/edges, the centroid and the set of edges of $E$, respectively. We further set, as usual, $h = \max_{E \in \Th} h_E$.

The virtual element space for the displacement is given by \cite{DaVeigaBrezzi2013}:
\begin{align*}
    \VMh[E]{1} = \Big\{\vv \in \vsob{1}{E}{2}: &(i)\ \Delta \vv = 0\\
    &(ii)\ \vv_{|e} \in \VPoly{1}{e}{2},\, \forall e \in \Eh[E] \text{ and } \vv \in [\con{0}{\partial E}]^{2}\Big\},    
\end{align*}
where $\Delta \vv = \begin{bmatrix}
    \Delta v_1\\
    \Delta v_2
\end{bmatrix}$ is the component-wise Laplacian operator.
The degrees of freedom for a vector-valued function $\vv \in \VMh[E]{1}$ are the values of $\vv$ at the vertices of the element $E$. Thus, we have
\begin{equation*}
    \Ndof[E] = 2 \Nv[E].
\end{equation*}

Let us define the set of Lagrange basis functions $\{\vvarphi_{i,E}\}_{i=1}^{\Ndof[E]}$ related to the aforementioned degrees of freedom. It is easy to check that
\begin{equation}
    \vvarphi_{i,E} = \begin{cases}
        \begin{bmatrix}
            \varphi_{i,E}\\
            0
        \end{bmatrix} & \text {if } i =1,\dots,\Nv[E],\\
        \begin{bmatrix}
            0 \\
            \varphi_{i - \Nv[E],E}\\
        \end{bmatrix} & \text {if } i =\Nv[E]+1,\dots,\Ndof[E],\\
    \end{cases}
    \label{eq:scalar_to_vec}
\end{equation}
where $\{\varphi_{j,E}\}_{j=1}^{\Nv[E]}$ corresponds to the set of Lagrange basis functions related to degrees of freedom for the scalar virtual element space \cite{LBe13}:
\begin{equation}
    \begin{aligned}
        \VPh[E]{1}{} = \Big\{v \in \sob{1}{E}: &(i)\ \Delta v = 0\\
        &(ii)\ v_{|e} \in \Poly{1}{e},\, \forall e \in \Eh[E] \text{ and } v \in \con{0}{\partial E} \Big\}. 
    \end{aligned}
    \label{eq:vem_scalr_space}
\end{equation}
Indeed, it holds $\VMh[E]{1} = [\VPh[E]{1}{}]^2$.

\subsection{The Virtual Element Discretization}\label{sec:vem_discretization}

For each polygon $E$, we introduce the \textit{computable} local polynomial projector $\proj{\nabla, E}{1} : \vsob{1}{E}{2} \to \VPoly{1}{E}{2}$, which is defined as
\begin{equation*}
\begin{cases}
    \scal[E]{\nabla \vv - \nabla \proj{\nabla, E}{1} \vv}{\nabla \pp} = 0\quad \forall \pp \in \VPoly{1}{E}{2}, \\
     \int_{\partial E} \proj{\nabla, E}{1} \vv = \int_{\partial E} \vv\quad \forall \vv \in \vsob{1}{E}{2}.
     \end{cases}
\end{equation*}
Moreover, we denote by $\proj{0,E}{0}$ both the vector-valued and the tensor-valued $L^2$-projectors onto constants. In particular, $\forall \vv \in \vsob{1}{E}{2}$, it holds \cite{DaVeigaBrezzi2023}
\begin{equation*}
    \proj{0,E}{0} \nabla \vv = \nabla \proj{\nabla,E}{1} \vv,
\end{equation*}
and 
\begin{equation*}
    \proj{0,E}{0} \div \vv = \proj{0,E}{0} \tr (\nabla \vv) = \tr (\proj{0,E}{0} \nabla \vv).
\end{equation*}

Specifically, by computable we mean that the projection of virtual functions can be computed only through the information provided by the degrees of freedom of such a function.
Now, let us split the continuous bilinear form \eqref{eq:continuous_a} according to $\Th$, i.e.
\begin{equation}
    \dbilin{\uu}{\vv} = \sum_{E \in \Th} \dbilin[E]{\uu}{\vv}, \quad \text{with}\quad \dbilin[E]{\uu}{\vv} = \int_E \ssigma(\xx, \nabla \uu) : \nabla \vv.
\end{equation}
For all $E \in \Th$, we define the discrete virtual bilinear form as \cite{DaVeiga2015}:
\begin{equation}
\begin{aligned}
    \nldbilinh[E]{\uu_h}{\vv_h}{\ww_h} &:= \int_E \ssigma(\xx, \proj{0,E}{0}\nabla \uu_h) : \proj{0,E}{0} \nabla \vv_h \\
    &\quad + \alpha_E(\ww_h) \stab[E]{(I-\proj{\nabla,E}{1})\uu_h}{(I-\proj{\nabla,E}{1})\vv_h}.
    \end{aligned}
    \label{eq:nl_discrete_bilin_form}
\end{equation}
where the \textit{stability term} $\stab[E]{}{}$ is chosen in such a way it scales like $\dbilin[E]{}{}$ for a unitary material constants on the kernel of $\proj{\nabla,E}{1}$. On the other hand, the \textit{stabilizing parameter} $\alpha_E(\ww_h)$ takes into account different material constants as well as non-linear materials.

Finally, to solve Problem \eqref{eq:nl_discrete_bilin_form}, we must address the issue of non-linearity. Among the commonly used numerical methods for solving systems of non-linear equations, the Newton–Raphson method is the fastest, offering a quadratic convergence provided that the initial estimate is sufficiently close to the solution. Therefore, selecting an initial estimate close to the solution is crucial for accelerating the convergence. In solid mechanics, the initial estimate is usually set to the undeformed shape of the structure, i.e. the initial displacement is set to the all-zeros vector. Since the initial estimate usually starts from the zero displacement, the Newton–Raphson method converges quickly when the applied load is small. However, convergence becomes challenging when large applied loads induce large displacements. The incremental force method addresses this issue by applying the load in increments. Within each increment, the Newton–Raphson procedure remains unchanged. A new load increment is applied after the solution corresponding to the previous load increment has successfully converged \cite{Kim2015}.

Thus, as in \cite{DaVeiga2015}, for each $n = 1, \dots, N$, we solve the following problem: \textit{Given $\ww_h \in \VMh{1}$, find $\uu_h^n \in \VMh{1}$ such that }
\begin{equation}
    \sum_{E \in \Th} \nldbilinh[E]{\uu_h^n}{\vv_h}{\ww_h} = \sum_{E \in \Th} \scal[E]{\ff^n}{\proj{0,E}{0}\vv_h} \quad \forall \vv_h \in \VMh{1},
    \label{eq:nl_discrete_var_problem}
\end{equation}
where $\ff^{n} = \frac{n}{N} \ff$ and $\uu_h^0 = \bm{0}$ on $\Omega$. For each $n$, Problem \eqref{eq:nl_discrete_var_problem} is solved with the Newton-Raphson method. Setting
\begin{equation*}
    \mathcal{G}_h^E(\uu_h^n, \vv_h;\ww_h) =  \nldbilinh[E]{\uu_h^n}{\vv_h}{\ww_h} -  \scal[E]{\ff^n}{\proj{0,E}{0}\vv_h},
\end{equation*}
we introduce 
\begin{equation}
    \begin{aligned}
    D \mathcal{G}_h(\uu_h^{n,k}, \vv_h;\ww_h) \cdot [\ddelta_h^{n,k}] &=   \int_{E} \left(\mathbb{A}(\xx, \proj{0,E}{0} \nabla \uu_h^{n,k}) : \nabla \proj{0,E}{0}\ddelta_h^{n,k} \right)  : \proj{0,E}{0}\nabla \vv_h \\
    &\quad+ \alpha_E(\ww_h) \stab[E]{(I-\proj{\nabla,E}{1})\ddelta_h^{n,k}}{(I-\proj{\nabla,E}{1})\vv_h},
\end{aligned}
\label{eq:linearization_vem}
\end{equation}
where the elastic modulus $\mathbb{A}(\xx, \nabla \uu)$ is defined as the fourth-order tensor
\begin{equation}
   \mathbb{A}(\xx, \nabla \uu) = \frac{\partial \ssigma(\xx, \nabla \uu)}{\partial \nabla \uu}.
   \label{eq:elastic_modulus}
\end{equation}
Finally, as in \cite{Wriggers2008}, we generate a sequence $\uu_h^{n,k+1} = \uu_h^{n,k} + \ddelta_h^{n,k}$ for each $k \geq 0$ by solving: \textit{Given $\ww_h \in \VMh{1}$, find $\ddelta_h^{n,k} \in \VMh{1}$ such that}
\begin{equation*}
   \sum_{E\in \Th} \left[\mathcal{G}^E_h(\uu_h^{n,k}, \vv_h;\ww_h) + D \mathcal{G}^E_h(\uu_h^{n,k}, \vv_h;\ww_h) \cdot [\ddelta_h^{n,k}] \right] = 0 \quad \forall \vv_h \in \VMh{1}.
\end{equation*}
Usually, $\ww_h$ is set equal to the solution found at the last loading increment step, i.e. $\uu^{n-1}_h$. This choice is taken to simplify the Newton-Raphson application by avoiding the computation of derivatives of $\alpha_E$ in the related tangent matrix \eqref{eq:linearization_vem}. Indeed, although $\ww_h = \uu_h^{n}$ is a more intuitive choice, it is also more computationally demanding \cite{DaVeiga2015, Chi2016}. Other possible choices include setting $\ww_h = \bm{0}$, i.e. evaluating $\alpha_E$ in the undeformed configuration. However, this might lead to unsatisfactory results, also in the small deformations regime \cite{Chi2016}.
Moreover, different choices for $\alpha_E$ have been introduced:
\begin{itemize}
    \item \textit{norm-based stabilization} \cite{DaVeiga2015}:
    \begin{equation*}
        \alpha_E(\ww_h) = \vertiii{\mathbb{A}(\xx, \proj{0,E}{0} \nabla \ww_h)},
    \end{equation*}
    which is chosen to produce a strictly positive stabilization term. The operator $\vertiii{\cdot}$ could be chosen as any fourth-order tensor norm. 
    \item \textit{trace-based stabilization} \cite{Chi2016}:
    \begin{equation*}
        \alpha_E(\ww_h) = \frac{1}{4} \tr\left(\mathbb{A}(\xx, \proj{0,E}{0} \nabla \ww_h)\right).
    \end{equation*}
    We observe that, with this choice, the stabilization may assume negative values.
    \item \textit{stiffness-based stabilization} \cite{Wriggers2024}:
    \begin{equation*}
        \alpha_E(\ww_h) = \frac{1}{4 \Ndof}  \left(\sum_{i=1}^{\Ndof} \left((\Kmatrix^C_T)_{ii}\right)^2\right)^{\frac{1}{2}} \quad \text{or}\quad \alpha_E(\ww_h) = \frac{1}{4 \Ndof}  \sum_{i=1}^{\Ndof} (\Kmatrix^C_T)_{ii},
    \end{equation*}
    where $\Kmatrix^C_T$ is the global tangent system matrix related to the consistency term.
\end{itemize}

\section{The Neural Approximated Virtual Element Method}\label{sec:navem}
In the NAVEM method, we propose to approximate $\varphi_{j,E}\in \VPh[E]{1}{}$ by means of a neural network-based approach and to assemble the final linear system as in the more classical Finite Element Method \cite{PintoreTeora2025}.  

To do so, we approximate $\varphi_{j,E}$ by a function $\nvarphi_{j,E}$ in a suitable functional space $\approxspace[j,E]$. A similar procedure is repeated to approximate the gradient of $\varphi_{j,E}$ (or, more generally, all the required derivatives). Finally, since the approximation $\nvarphi_{j,E}$ of $\varphi_{j,E}$ is known in closed form as well as its derivatives, it can be evaluated point-wise inside each element $E\in\Th$ to numerically compute the integrals in \eqref{eq:var_problem}, without introducing any projection or stability operator.

\subsection{The Local Approximation Space}\label{sec:local_space}

In this section, we describe the local function space $\approxspace[j,E]$ in which the approximation $\nvarphi_{j,E}$ of $\varphi_{j,E}$ is sought. Such space coincides with the one adopted in \cite{PintoreTeora2025}. Here, for the sake of clarity, we briefly recall the main concepts and refer to \cite{PintoreTeora2025} for a more detailed description.

All the functions in $\approxspace[j,E]$ are defined on a reference squared region $\refregion = [-\refedge, \refedge]^2\subset \R^2$ independent of $E$. Therefore, we assume that a suitable affine map is available to map each element $E\in\Th$ to a reference element $\hat E \subset \refregion$. A convenient map is the one described in \cite{Teora2024}, but other choices are admissible.

Let $\HPoly{\refell}{\refregion}$ be the space of harmonic polynomials up to order $\refell\ge0$. A basis of such a space is the set of scaled polynomials
\begin{equation}
      \Big\{ 1,\ \Re\left(\left(\frac{z}{\refedge}\right)^\ell\right), \Im\left(\left(\frac{z}{\refedge}\right)^\ell\right),\ \ell = 1,\dots, \refell \Big\}.
      \label{eq:scaled_harmonic_polynomial}
\end{equation}

Here $z$ is the complex number $z=x_1+ix_2$ associated with the point $\xx=[x_1,x_2]^T\in\R^2$, and $\Re$ and $\Im$ denote the real and imaginary part of a complex quantity. Note that all the polynomials in \eqref{eq:scaled_harmonic_polynomial} can be recursively retrieved as in \cite{Perot2021}.

As discussed in \cite{PintoreTeora2025}, harmonic polynomials are not enough to accurately approximate the VEM basis functions, we thus introduce additional functions mimicking the target virtual functions. Let us introduce the Laplace problem
\begin{equation*}
\begin{cases}
        \Delta \tilde{\Phi} = 0 & \text{in } \Omega_{\Phi} = (-1,1)^2,\\
        \tilde{\Phi} = 1 + x_2 & \text{on } \Gamma_{\Phi, 1} = \{x_1 = 1 \text{ and } -1 \leq x_2 \leq 0\},\\
        \tilde{\Phi} = 1 - x_2 & \text{on } \Gamma_{\Phi, 2} = \{x_1 = 1 \text{ and } 0 \leq x_2 \leq 1\},\\
        \tilde{\Phi} = 0 & \text{on } \partial \Omega_{\Phi} \setminus \{\Gamma_{\Phi, 1} \cup \Gamma_{\Phi, 2}\},
\end{cases}
\end{equation*}
and approximate its solution with a function $\Phi$ of the form
\begin{equation}
    \Phi(z) = \sum_{\alpha =1}^{N^{\text{1}}} c^{\text{1}}_{\alpha} \Re\left(\frac{d_{\alpha}}{z - z_{\alpha}}\right) + \sum_{\beta =0}^{N^{\text{2}}} c^{\text{2}}_{\beta} \Re\left(\left(\frac{z}{2}\right)^{\beta}\right).
    \label{eq:hanging_function}
\end{equation}
In this formula, $\Re\left(\left(\frac{z}{2}\right)^{\beta}\right)$, $\beta=0,\dots,N^2$, are known harmonic polynomials and $\Re\left(\frac{d_{\alpha}}{z - z_{\alpha}}\right)$ are known harmonic rational functions associated with known poles $z_\alpha = 1 + 2 \exp\left(-4 (\sqrt{N_1} - \sqrt{\alpha})\right)$, for $\alpha=1,\dots,N^{\text{1}}$. The sought coefficients $c_\alpha^1$, $\alpha=1,\dots,N^1$, and $c_\beta^1$, $\beta=1,\dots,N^2$, are computed by solving a linear least square problem to minimize the difference between $\Phi$ and $\tilde \Phi$ on a given set of points on $\partial\Omega_\Phi$. The function $\Phi$ is then used to compute three auxiliary functions $\Phi_{j,E}^{j-1}$, $\Phi_{j,E}^{j}$, and $\Phi_{j,E}^{j+1}$ through three different affine maps. Such maps map the point $(1,0)$ to the $(j-1)$, $j$, or $(j-1)$-th vertex of $E$, and ensure that $E$ is included in the image of $\Omega_\Phi$. We note that, since $\Phi$ does not depend on $E$ or $\widehat{E}$, it is computed only once.

Finally, for any element $E\in\Th$ and each vertex $v_j$ of $E$, with $=1,\dots,\Nv[E]$, we can define the space $\approxspace[j,E]$,  as:
\begin{equation}
    \approxspace[j,E] = \myspan \Big\{ \{\tilde{p}_{\beta}\}_{\beta=1}^{2\refell + 1},\ \Phi^{j-1}_{j,E}, \Phi^{j}_{j,E}, \Phi^{j+1}_{j,E} \Big\},
    \label{eq:approxspace}
\end{equation}
where $\{\tilde{p}_{\beta}\}_{\beta=1}^{2\refell + 1}$ are the harmonic polynomials in \eqref{eq:scaled_harmonic_polynomial}. Therefore, the dimension of $\approxspace[j,E]$ is
\begin{equation}
    \dim \approxspace[j,E] = 2\refell + 4.
\end{equation}
We observe that this dimension does not depend on the number of functions $N^1+N^2$ used to approximate $\Phi$.

\subsection{The Neural Network Architecture}\label{sec:nn_architecture}

As previously stated, the main ingredient of the NAVEM method is the approximation $\nvarphi_{j,E}$ of the local VEM basis function $\varphi_{j,E}$. This approximation is a suitable function in $\approxspace[j,E]$ which is built through a neural network trained by mimicking the function $\varphi_{j,E}$ at the boundary of $E$, where such function is well-known. For the ease of notation, we define
\begin{equation}\label{eq:approxspace_h}
\approxspace[j,E] = \myspan\{\mathit{h}_k\}_{k=1}^{\dim \approxspace[j,E]},
\end{equation}
where $\{\mathit{h}_k\}_{k=1}^{\dim \approxspace[j,E]}$ coincides with the set of basis functions in \eqref{eq:approxspace}. Therefore, the function $\nvarphi_{j,E}$ can be expressed as
\begin{equation}\label{eq:nvarphi_lin_comb}
\nvarphi_{j,E} = \sum_{k=1}^{\dim \approxspace[j,E]} c_k^\varphi(v_j,E) h_k
\end{equation}
for suitable coefficients $c_k^\varphi(v_j,E)$ depending on the pair $(v_j,E)$. In this context, the involved neural network is a parametric function mapping the pair $(v_j,E)$ to the vector $[c_k^\varphi(v_j,E)]_{k=1}^{\dim \approxspace[j,E]}$. We remark that, since a neural network accepts only vectors as inputs, the pair $(v_j,E)$ has to be suitably encoded. Multiple input encoding strategies are admissible. In this work, we adopt the approach described in \cite{PintoreTeora2025}, where an input encoding strategy is detailed, along with variability and input reduction strategies for generic polygons.

Let $\xx_0$ be the encoding of a pair $(v_j,E)$. In the following, we drop the subscripts $j,E$ to lighten the notation, while also admitting that the neural network can be evaluated on multiple pairs at once. The operator mapping $\xx_0$ to the coefficients ${\mathbf c^{\varphi}}(\xx_0):= [c_k^\varphi]_{k=1}^{\dim \approxspace[j,E]}$ is expressed as:
\begin{equation} 
\begin{aligned}
 & \xx_0^{\varphi} := \xx_0, \\
 &\xx_\ell^{\varphi} = \rho(A^{\varphi}_\ell \xx_{\ell-1}^{\varphi} + b^{\varphi}_\ell), \hspace{2cm} \ell = 1,...,L-1, \\
 &{\mathbf c^{\varphi}}(\xx_0) = A^{\varphi}_{L} \xx_{L-1}^{\varphi}  + b^{\varphi}_L.
  \end{aligned}
  \label{eq:nn_formula}
\end{equation}

Equation \eqref{eq:nn_formula} describes the so-called architecture of the neural network \cite{goodfellow2016deep}. Here, $L$ denotes the number of layers, $\rho$ is a non-linear activation function, and $A^{\varphi}_\ell\in\R^{N_\ell\times N_{\ell-1}}$ and $b^{\varphi}_\ell\in\R^{N_\ell}$ are matrices and vectors containing the trainable coefficients of the neural network, which are optimized by minimizing a suitable loss function. Specifically, to optimize the coefficients of this neural network, we propose to minimize the distance between the traces of the functions $\nvarphi_{j,E}$ and $\varphi_{j,E}$ on $\partial E$, i.e. we want to minimize the following quantity
\begin{equation}
\epsilon_{j,E}^\varphi = \norm[{\sob{1/2}{\partial E}}]{ \nvarphi_{j,E} - \varphi_{j,E}},
\label{eq:approximation_error}  
\end{equation}
for each pair $(j,E)$ in a given training dataset. We remark that this quantity is computable because it involves only the evaluations of $\varphi_{j,E}$ on the boundary of $E$, where this function is known in closed form.

Let $\qq_{j,E}:=\nabla \varphi_{j,E}$ be the gradient of $\varphi_{j,E}$. We observe that his gradient can be approximated by the gradient of the function $\nvarphi_{j,E}$ defined in  \eqref{eq:nvarphi_lin_comb} \cite{PintoreTeora2024}, which can be exactly computed as:
\begin{equation}\label{eq:derivative_lin_combin_v0}
\nabla \nvarphi_{j,E} = \sum_{k=1}^{\dim \approxspace[j,E]} c_k^\varphi(v_j,E) \nabla h_k.
\end{equation}
However, such a strategy is observed to be suboptimal. For this reason, it is convenient to introduce a second neural network to approximate $\qq_{j,E}$ with vector-functions contained in $\nabla \approxspace[j,E]$. This neural network, similarly to \eqref{eq:nn_formula}, can be expressed as:
\begin{equation} 
\begin{aligned}
 & \xx_0^{\qq} := \xx_0, \\
 &\xx_\ell^{\qq} = \rho(A^{\qq}_\ell \xx_{\ell-1}^{\qq} + b^{\qq}_\ell), \hspace{2cm} \ell = 1,...,L-1, \\
 &{\mathbf c^{\qq}}(\xx_0) = A^{\qq}_{L} \xx_{L-1}^{\qq}  + b^{\qq}_L.
  \end{aligned}
  \label{eq:nn_formula_der}
\end{equation}
The output ${\mathbf c^{\qq}}(\xx_0):= [c_k^{\qq}]_{k=1}^{\dim \approxspace[j,E]}$ of \eqref{eq:nn_formula_der} represents the set of coefficients that allows us to approximate $\qq_{j,E}$ through the following vector 
\begin{equation}\label{eq:derivative_lin_combin_v2}
\qq_{j,E}^\NN = \sum_{k=1}^{\dim \approxspace[j,E]} c_k^{\qq}(j,E) \nabla h_k.
\end{equation}
The trainable coefficients $A^{\qq}_{\ell}$ and $b^{\qq}_{\ell}$ are optimized by minimizing a \textit{computable} loss function. This cost term is defined through the following quantity
\begin{equation}
\epsilon_{j,E}^{\qq} = \norm[{\leb{2}{\partial E}}]{ \left(\qq_{j,E}^\NN  - \qq_{j,E}\right) \cdot {\bm{t}}},
\label{eq:approximation_error_der}  
\end{equation}
where $\bm{t}$ denotes the tangential vector to the boundary $\partial E$. In \cite{PintoreTeora2025}, the authors proved that small values of $\epsilon_{j,E}^{\varphi}$ and $\epsilon_{j,E}^{\qq}$ are enough to ensure a good approximation of $\varphi_{j,E}$ and of $\qq_{j,E}$ on the entire element $E$, even if they are minimized only on $\partial E$.
We observe that such an approach introduces a consistency error because \eqref{eq:derivative_lin_combin_v2} does not exactly coincide with the derivative of \eqref{eq:nvarphi_lin_comb}, but it significantly reduces the oscillations in the gradients of the NAVEM basis functions, leading to a more stable method.

\subsection{The Training Phase}\label{sec:nn_training}

We highlight that employing the input encoding strategy presented in \cite{PintoreTeora2025} requires defining a different network for each class of polygons characterized by a different number of vertices $\Nv \geq 4$. For $\Nv = 3$, the virtual element basis functions, as well as the NAVEM basis functions, coincide with the ones of the standard FEM. 
Thus, each training set has to contain only the class of polygons corresponding to the network that is trained. 

For each class of polygons, let $\mathcal{T}^{\mathrm{train}}$ be a training dataset of size $\#  \mathcal{T}^{\mathrm{train}}$ containing polygons $E$ with $\Nv$ vertices. Then, to train the neural networks \eqref{eq:nn_formula} and \eqref{eq:nn_formula_der}, i.e. to approximate the basis functions and their gradients, we minimize the following loss functions 
\begin{equation}
\begin{gathered}
    {\cal L}_\varphi = \sqrt{\frac{1}{\#  \mathcal{T}^{\mathrm{train}} \Nv} \sum_{E\in {\mathcal{T}^{\mathrm{train}}}} \sum_{j=1}^{\Nv} \left( \epsilon_{j,E}^\varphi \right)^2}, \\
    {\cal L}_{\qq} = \sqrt{\frac{1}{\#  \mathcal{T}^{\mathrm{train}} \Nv} \sum_{E\in {\mathcal{T}^{\mathrm{train}}}} \sum_{j=1}^{\Nv} \left( \epsilon_{j,E}^{\qq} \right)^2},
    \end{gathered}
    \label{eq:used_losses}
\end{equation}
respectively. Finally, a standard regularization term penalizing the $L^2$-norm of the trainable weights, with regularization coefficients $10^{-8}$, is added to both losses. We observe that these formulas differ from the ones defined in  \cite{PintoreTeora2025}. Indeed, in \cite{PintoreTeora2025}, the used loss functions coincide with the squared losses ${\cal L}_\varphi^2$ and ${\cal L}_{\qq}^2$. Note that the global minima of these two pairs of losses are the same, but we empirically observed better performance minimizing \eqref{eq:used_losses}. Moreover, we here decide to employ a new optimizer: after performing $5000$ epochs with the first-order ADAM optimizer \cite{kingma2014adam} as in \cite{PintoreTeora2025}, we carry out $5000$ updates with the self-scaled BFGS optimizer described in \cite{wright1999numerical, kiyani2025optimizer}. This last advanced optimization technique enhances convergence by adaptively scaling the weights updates provided by the BFGS optimizer. 

To appreciate the effectiveness of these new techniques, we compare the final values of the training losses by training the neural networks as described here and in the previous manuscript \cite{PintoreTeora2025}. More precisely, we here build new neural networks characterized by the same architecture used in \cite{PintoreTeora2025}, i.e. they all have 5 layers, with 50 neurons in each layer, and the \textit{tanh} activation function in each hidden layer. For the network associated with the class of polygons with $\Nv$ vertices, the input dimension is $2(\Nv-1)$, while the output is always $2\refell+4$, as discussed in Section \ref{sec:local_space}. We fix $\refell=20$, independently of the associated class of polygons, and set $\refedge = 3$. A fine-tuning strategy of all these values should be advisable, but it is beyond the scope of the current paper. Moreover, in order to train such networks, we use the same datasets RDQM, for $\Nv =4$, and VM, for $\Nv \geq 5$, presented in \cite{PintoreTeora2025}: the former contains randomly generated quadrilaterals, and the latter contains polygons sampled from Voronoi meshes generated through \cite{mvem}. The entire code is written in Python, and the implementation of the definition, optimization, and evaluation of the neural networks heavily relies on the TensorFlow library \cite{tensorflow2015-whitepaper}. 

The differences between the accuracy of the two approaches can be appreciated in Table \ref{tab:losses_comparison}. In particular, we report the value of $\mathcal{L}_{\varphi}$ and $\mathcal{L}_{\qq}$ for the networks used in \cite{PintoreTeora2025} and for the networks used in the current manuscript. It is evident that the change of the loss and of the optimizer allows one to more accurately approximate both the basis functions and their gradients. We remark that these changes do not influence the overall computational cost of the training phase because the cost of the additional operations is negligible.

\begin{table}[!ht]
\centering
\caption{Final training loss functions computed using the neural networks in \cite{PintoreTeora2025} and the ones in the current manuscript.}
\label{tab:losses_comparison}
\resizebox{0.5\textwidth}{!}{%
\begin{tabular}{@{}ccccc@{}}
\toprule
    \multicolumn{1}{c}{$\Nv$}              & $\mathcal{L}_{\varphi}$ \cite{PintoreTeora2025} & $\mathcal{L}_{\varphi}$       & $\mathcal{L}_{\qq}$ \cite{PintoreTeora2025} & $\mathcal{L}_{\qq}$ \\ \midrule
\multicolumn{1}{c|}{4} & 4.62e-03                                        & \multicolumn{1}{c|}{5.13e-04} & 1.00e-02                                    & 2.94e-03            \\
\multicolumn{1}{c|}{5} & 1.97e-03                                        & \multicolumn{1}{c|}{2.81e-04} & 4.97e-03                                    & 1.84e-03            \\
\multicolumn{1}{c|}{6} & 1.02e-03                                        & \multicolumn{1}{c|}{1.12e-04} & 2.96e-03                                    & 1.07e-03            \\
\multicolumn{1}{c|}{7} & 2.15e-03                                        & \multicolumn{1}{c|}{3.40e-04} & 5.22e-03                                    & 2.07e-03 \\ 
\multicolumn{1}{c|}{8} & 1.21e-03                                        & \multicolumn{1}{c|}{3.26e-04} & 4.88e-03                                    & 2.47e-03  \\ \bottomrule
\end{tabular}%
}
\end{table}

We would like to remark that the comparison is meaningful because we are aiming to approximate the same scalar virtual element space $\VPh[E]{1}{}$.

\subsection{The Online Phase and the NAVEM Discretization}\label{sec:navem_scheme}

\begin{figure}[!ht]
\centering
\includegraphics[width=\textwidth]{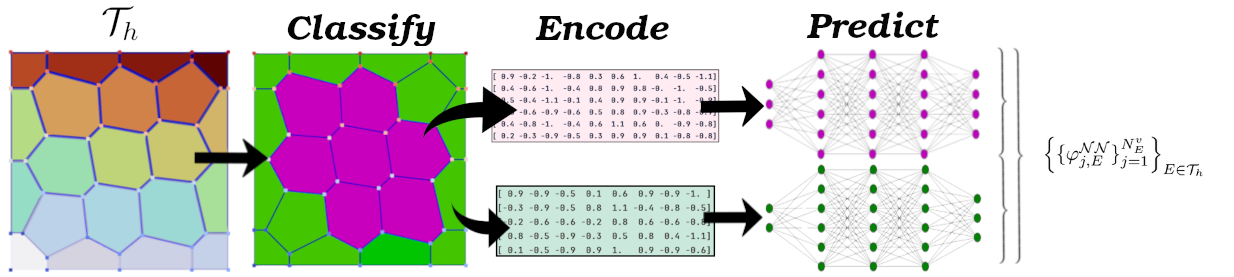}
\caption{A sketch of the NAVEM online phase (steps \ref{step:classify}--\ref{step:predict}).}
\label{fig:online_pahse}
\end{figure}

In this section, assuming that suitable trained neural networks are available to approximate the basis functions related to $\Th$, we present the NAVEM discretization of the elastic problem \eqref{eq:model_problem}. The main steps of the entire procedure can be subdivided in two main parts. The first one (steps \ref{step:classify}--\ref{step:predict}) is specific to the NAVEM approach and it is summarized in Figure \ref{fig:online_pahse}. The second one (steps \ref{step:assemble}--\ref{step:solve}) depends only on the problem and can be handled as in the standard Finite Element Method.

Let us describe the first part. For each element $E$ in the tessellation $\Th$ of the domain $\Omega$, and for each vertex $v_j$ of $E$, with $j = 1,\dots,\Nv[E]$, we
\begin{enumerate}[label=\textbf{S.\arabic*}]
\item \label{step:classify}\textit{classify} the pair $(v_j, E)$, as discussed in Section \ref{sec:nn_training}, to select the neural networks that have to be used to approximate the corresponding basis function and its gradient; 
\item \textit{encode} the information $(v_j,E)$ to define a vector $\xx_0^{j,E}$ that represents the input vector for the neural networks;
\item \label{step:predict} \textit{predict} the coefficients ${\mathbf c^{\varphi}}(v_j, E)$ and ${\mathbf c^{\qq}}(v_j, E)$ to compute $\nvarphi_{j,E}$ and $\qq_{j,E}^\NN$ according to \eqref{eq:nvarphi_lin_comb} and \eqref{eq:derivative_lin_combin_v2}, respectively.
\end{enumerate}

Now, we define the local scalar NAVEM space $\nVh[E]{1}$ and the related gradient space $\nabla \nVh[E]{1}$ as:
\begin{equation*}
\begin{gathered}
    \nVh[E]{1} = \myspan\left\{\nvarphi_{j,E}, \,\,j=1,\dots,\Nv[E]\right\},\\
    \nabla \nVh[E]{1} = \myspan\left\{\qq^{\NN}_{j,E}, \,\,j=1,\dots,\Nv[E]\right\},
\end{gathered}
\end{equation*}
which represent the approximation spaces for the scalar VEM space $\VPh[E]{1}{}$, defined in \ref{eq:vem_scalr_space}, and the related gradient space $\nabla \VPh[E]{1}{}$, respectively. 

The local NAVEM space for the displacement and its gradient can be defined following a standard strategy, i.e. it can be defined as
\begin{equation}
    \nVMh{1}(E) := \nVh[E]{1} \times \nVh[E]{1},\quad \nabla \nVMh[E]{1} = \nabla \nVh[E]{1} \times \nabla \nVh[E]{1}.
\end{equation}
for each element $E \in \Th$. By gluing together these local spaces, we obtain the global NAVEM discrete spaces for the displacement. Finally, we observe that, even if the NAVEM functions are not exactly continuous across elements, we do not double the degrees of freedom. Thus,
\begin{equation}
    \Ndof = \dim(\nVMh{1})=\dim(\VMh{1}),
\end{equation}
i.e. the dimension of $\nVMh{1}$ is twice the number of vertices of $\Th$ that do not belong to the Dirichlet boundary. 

Specifically, each NAVEM function $\vv_h^{\NN} \in \nVMh {1}$ and its gradient can be written as:
\begin{gather*}
    \vv_h^{\NN} = \sum_{i=1}^{\Ndof/2} v_i \begin{bmatrix}
        \nvarphi_i\\
        0
    \end{bmatrix} + \sum_{i=\Ndof/2}^{\Ndof} v_i \begin{bmatrix}
        0\\
        \nvarphi_{i - \Ndof/2}
    \end{bmatrix},\\ 
    \nabla \vv_h^{\NN} = \sum_{i=1}^{\Ndof/2} v_i \begin{bmatrix}
        \qq^{\NN}_{1,i} & \qq^{\NN}_{2,i}\\
        0 & 0
    \end{bmatrix} + \sum_{i=\Ndof/2}^{\Ndof} v_i \begin{bmatrix}
            0 & 0\\
        \qq^{\NN}_{1,i - \Ndof/2} & \qq^{\NN}_{2,i - \Ndof/2}
    \end{bmatrix},
\end{gather*}
where $\{v_i\}_{i=1}^{\Ndof}$ represents the set of degrees of freedom for the vector-function $\vv_h^{\NN}$.

Given the output of the online phase \ref{step:classify}--\ref{step:predict}, the NAVEM implementation follows the exact steps performed to solved the elastic problem \eqref{eq:var_problem} with a standard Finite Element Method. In particular, the NAVEM discretization for the elastic problem \eqref{eq:var_problem} reads as: \textit{Find $\uu_h^{\NN} \in \nVMh{1}$ such that}
\begin{equation*}
    \sum_{E \in \Th} \scal[E]{\ssigma(\xx,\nabla \uu_h^{\NN})}{\nabla \vv_h^{\NN}} = \sum_{E \in \Th} \scal[E]{\ff}{\vv_h^{\NN}} \quad \forall \vv_h^{\NN} \in \nVMh{1}.
\end{equation*}
This problem is a non-linear problem that can be solved using a Newton-Raphson solver, eventually paired with an incremental force strategy when large displacements are taken into account. Thus, given an initial guess $\uu^{\NN,0} \in \nVMh{1}$, for each non-linear step $k \geq 0$, we must solve the following problem: \textit{Find $\ddelta_h^{\NN, k} \in \VMh{1}$ such that}
\begin{equation*}
   \sum_{E\in \Th} \left[\mathcal{G}^E_{\NN}(\uu_h^{\NN, k}, \vv^{\NN}_h) + D \mathcal{G}^E_{\NN}(\uu_h^{\NN, k}, \vv_h^{\NN}) \cdot [\ddelta_h^{\NN, k}] \right] = 0 \quad \forall \vv_h \in \VMh{1},
\end{equation*}
where
\begin{equation}
    \mathcal{G}^E_{\NN}(\uu_h^{\NN, k}, \vv^{\NN}_h) = \int_{E} \ssigma(\xx, \nabla \uu_h^{\NN, k}) : \vv_h^{\NN} - \int_E \ff \cdot \vv_h^{\NN},
    \label{eq:nn_rh}
\end{equation}
\begin{equation}
    D \mathcal{G}^E_{\NN}(\uu_h^{\NN, k}, \vv_h^{\NN}) \cdot [\ddelta_h^{\NN, k}] =  \int_{E} \left( \mathbb{A}(\xx, \nabla \uu_h^{\NN, k}) : \nabla \bm{\delta}^{\NN,k}_h \right): \nabla \vv_h^{\NN},
    \label{eq:nn_tangent_stiff}
\end{equation}
and $\mathbb{A}$ is the elastic modulus defined in \eqref{eq:elastic_modulus}. In the following, we will denote with the same bold symbol both the vector-function $\vv_h^{\NN} \in \nVMh{1}$ and the vector of the related degrees of freedom, assuming the meaning is clear from the context. 
In conclusion, the second part of the NAVEM discretization can be summarized as:
\begin{enumerate}[label=\textbf{S.\arabic*}]
\setcounter{enumi}{3}
\item\label{step:assemble} \textit{assemble} the tangent stiffness matrix $\mathbf{K}^{\NN, k}$ corresponding to \eqref{eq:nn_tangent_stiff} and the right-hand side $\mathbf{r}^{\NN, k}$ related to \eqref{eq:nn_rh}.
\item \label{step:solve} \textit{solve} the system 
\begin{equation}
    \mathbf{K}^{\NN, k} \ddelta^{\NN,k} = -\mathbf{r}^{\NN, k},
\end{equation}
and set $\uu^{\NN,k+1} = \ddelta^{\NN,k} + \uu^{\NN,k}$.
\end{enumerate}
The process is repeated until convergence. The exit vector $\uu^{\NN,k+1}$ collects the coefficients $\{u_i\}_{i=1}^{\Ndof}$ that describe the NAVEM solution $\uu^{\NN}_h$.


\section{Numerical Results}\label{sec:numerical_results}

In this section, we test the discussed method on linear and non-linear elasticity problems and on different meshes. 

We highlight that in all the experiments, the NAVEM basis functions are predicted by exploiting the same neural networks built to perform the experiment in Section \ref{sec:nn_training}. 

Denoting by $\uu$ the exact displacement, we test the performance of the NAVEM by looking at the behaviour of the following errors
\begin{equation}
\begin{gathered}
    \mathrm{err}_0^\NN = \sqrt{\sum_{E\in\Th} \norm[\leb{2}{E}]{\uu - \uu^{\NN}_h}^2},\\ \mathrm{err}_1^\NN = \sqrt{\sum_{E\in\Th} \norm[\leb{2}{E}]{\nabla \uu - \nabla \uu^{\NN}_h}^2},
    \end{gathered}
    \label{eq:nn_errors}
\end{equation}
as the mesh parameter $h$ decreases. For comparison purposes, we solve these problems also with the standard virtual element method. In this case, since we can not access the point-wise evaluation of virtual functions, we define the VEM errors as usual, that is \cite{LBe16}
\begin{equation}
\begin{gathered}
    \mathrm{err}_0^{\rm{VEM}} = \sqrt{\sum_{E\in\Th} \norm[\leb{2}{E}]{\uu - \proj{E,0}{1} \uu_h^{\rm{VEM}}}^2},\\
    \mathrm{err}_1^{\rm{VEM}} = \sqrt{\sum_{E\in\Th} \norm[\leb{2}{E}]{\nabla \uu - \proj{E,0}{0} \nabla \uu_h^{\rm{VEM}}}^2},
        \end{gathered}
    \label{eq:vem_errors}
\end{equation}
where we denote by $\uu_h^{\rm{VEM}}$ the virtual element solution.
We remark that to evaluate the high-order projection $\proj{E,0}{1} \uu_h^{\rm{VEM}}$ in virtual element, an enhanced version of the local space should be taken into account \cite{Ahmad2013}.

\subsection{Test 1: The linear case - convergence test}

\begin{figure}[!ht]
    \centering
    \begin{subfigure}{0.3\linewidth}
    \includegraphics[width=1\linewidth]{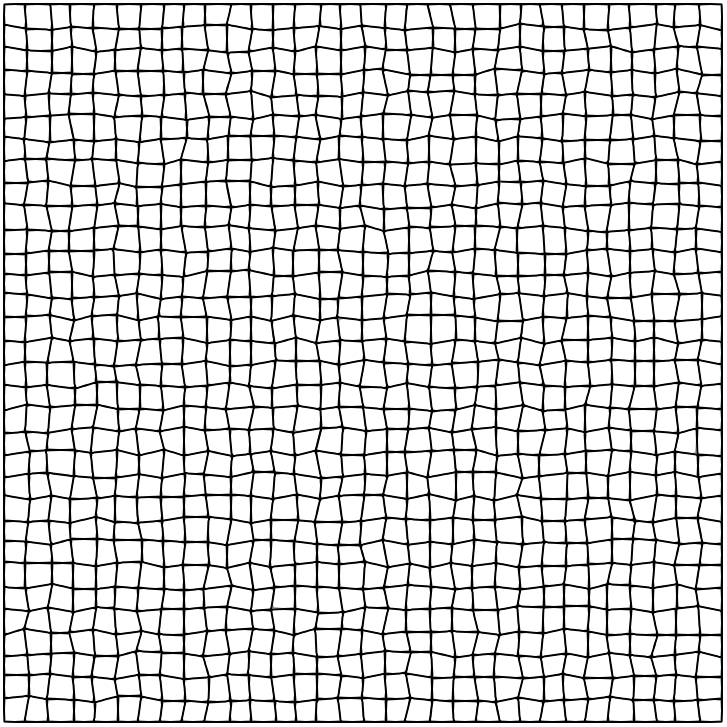}
    \caption{}
    \end{subfigure} \hspace{30pt}   
    \begin{subfigure}{0.35\linewidth}
    \includegraphics[width=1\linewidth]{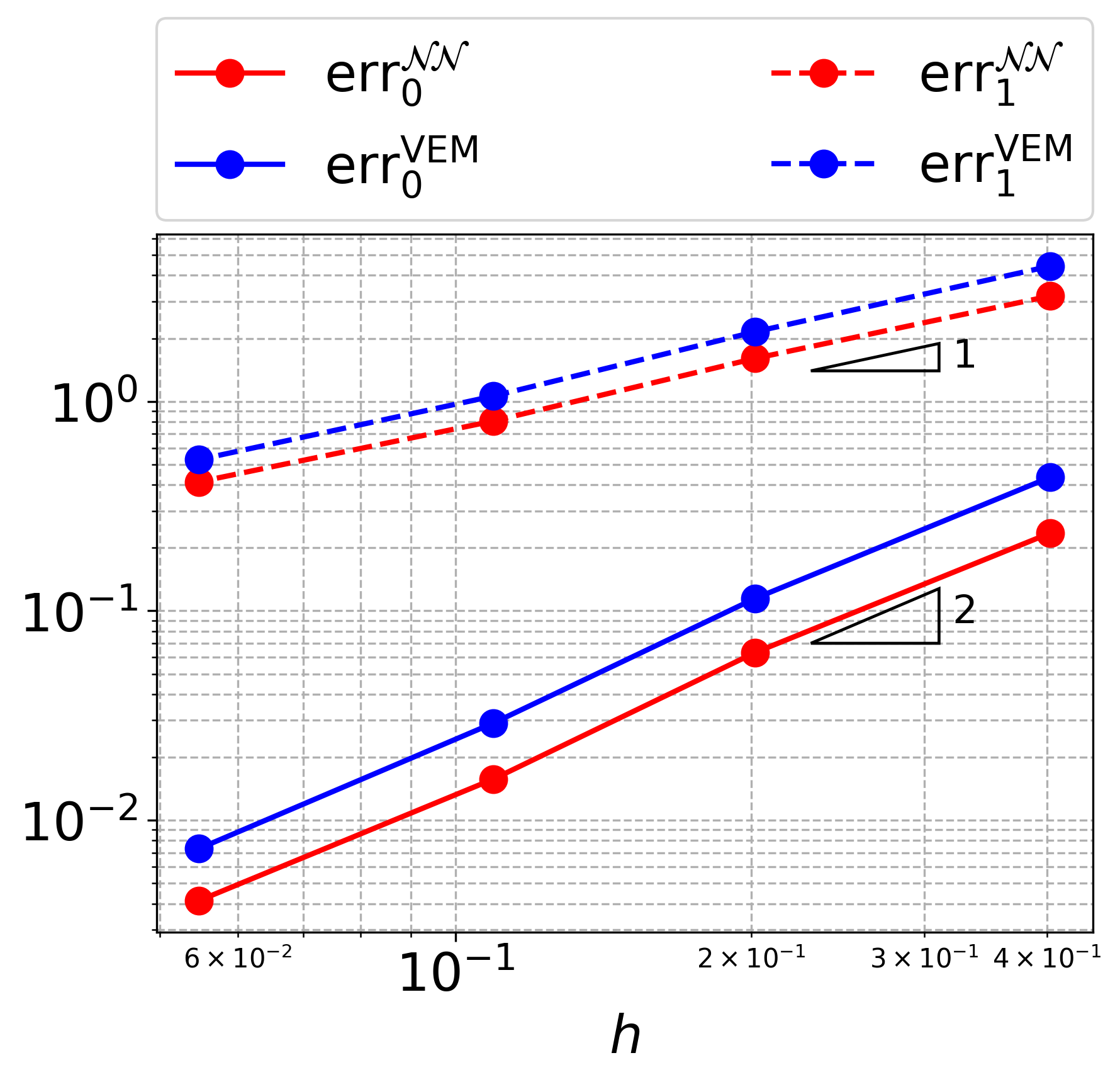}
    \caption{}
    \end{subfigure}
    \caption{Test 1: Left: Last refinement of the Distorted-Square family of meshes. Right: Behaviour of errors \eqref{eq:nn_errors}--\eqref{eq:vem_errors} as the mesh parameter $h$ decreases.}
    \label{fig:rnd_distorte_03_32x32}
\end{figure}

We consider the square domain $\Omega=(0,1)^2$ and the linear version of the problem \eqref{eq:var_problem}, where 
\begin{equation}
    \dbilin{\uu}{\vv} = 2 \mu \int_{\Omega} \eepsilon(\uu) : \eepsilon(\vv) + \lambda \int_{\Omega} \div \uu\ \div \vv,
\end{equation}
and $\mu, \lambda$ are the positive Lamé coefficients, whereas $\eepsilon(\uu) = \frac{1}{2}(\nabla \uu + (\nabla \uu)^T)$. We highlight that the algorithm detailed in the previous sections is valid for any constitutive law and, therefore, it offers a flexible setting to deal with elastic problems \cite{DaVeiga2015}.

In this test, we take the material constant values $\mu = 1.5$, $\lambda = 3$ and we set the boundary conditions and the body load $\ff$ in such a way the exact displacement is:
\begin{equation}
    \uu(\xx) = \begin{bmatrix}
        u_1(\xx)\\
        5 u_1(\xx)
    \end{bmatrix}, \qquad u_1(\xx) = 16 x_1(1-x_1)x_2(1-x_2) + 1.1.
\end{equation}
In particular, we choose to clamp the displacement at $\Gamma_D = \{\xx \in \Gamma: x_2 = 0\} \cup \{\xx \in \Gamma: x_1 = 0\}$, whereas the remaining boundary is free.

We build a family of four distorted quadrilateral meshes obtained by randomly distorting square grids of $4 \times 4$, $8 \times 8$, $16 \times 16$, and $32 \times 32$ elements. The last refinement is shown in Figure \ref{fig:rnd_distorte_03_32x32}. In the Virtual Element Method, the stabilizing parameter is set to be $2 \mu$. 

The behavior of errors \eqref{eq:nn_errors}-\eqref{eq:vem_errors} are shown in Figure \ref{fig:rnd_distorte_03_32x32}. We observe that the expected order of convergence for both methods is achieved, i.e. the errors decrease as $O(h^2)$ for the $L^2$-error and as $O(h)$ for the error in the $H^1$-seminorm.
Moreover, coherently with the results obtained for the scalar version of the NAVEM method in \cite{PintoreTeora2025}, we observe lower error constants with respect to VEM for both errors. We also note that the elements in this family of meshes are well represented by the training set. Thus, the network accuracy does not influence these results. 

Finally, we want to remark that the evaluation of a neural network has an intrinsic cost that tends to be negligible for fine meshes, where the cardinality of the input of the neural networks grows, as shown in \cite{PintoreTeora2025}. In a vectorial setting, as the one offered by the elasticity problems, this intrinsic cost is even lower because the neural networks are only used to construct the basis functions of $\nVh[E]{1}$. When a function in $\VM^\NN_{h,1}(E)$ has to be evaluated, the functions in $\nVh[E]{1}$ associated with each dimension are reused multiple times. This is possible because the vector space is just the Cartesian product of a scalar space, approximated with the neural network, with itself.

\subsection{Test 2: A benchmark non-linear elastic problem with analytical solution}

In this section, we propose a test taken from \cite{DaVeiga2015} to show some advantages of using an approach that does not require any stabilization term, such as the NAVEM. Specifically, we consider $\Omega = (0,1)^2$ and set $\Gamma_D = \Gamma$. We select the following constitutive law
\begin{equation}
    \ssigma(\xx, \nabla \uu (\xx)) = \mu(\bm{\epsilon}(\uu)) \bm{\epsilon}(\uu) 
\end{equation}
where the scalar function $\mu$ is given by
\begin{equation}
    \mu(\bm{\epsilon}(\uu)) = 3 (1 + \norm{\bm{\epsilon}(\uu)}^2) = 3 (1 + \sum_{n,m=1}^2 [\bm{\epsilon}(\uu)_{nm}]^2).
\end{equation}

\begin{figure}[!h]
    \centering
    \begin{subfigure}{0.3\textwidth}
    \includegraphics[width=1\textwidth]{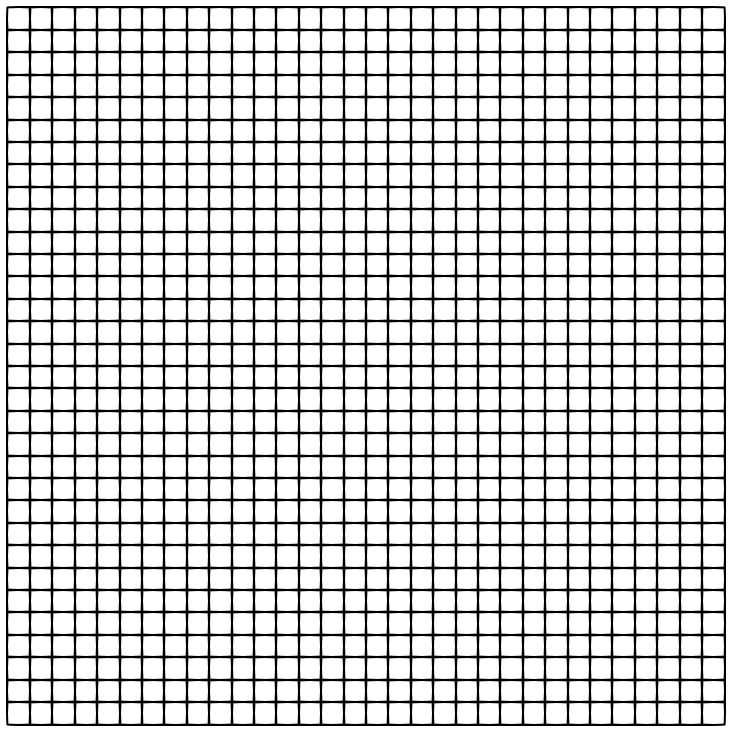}
    \caption{}
    \end{subfigure}    \hspace{30pt}  
    \begin{subfigure}{0.3\textwidth}
    \includegraphics[width=1\textwidth]{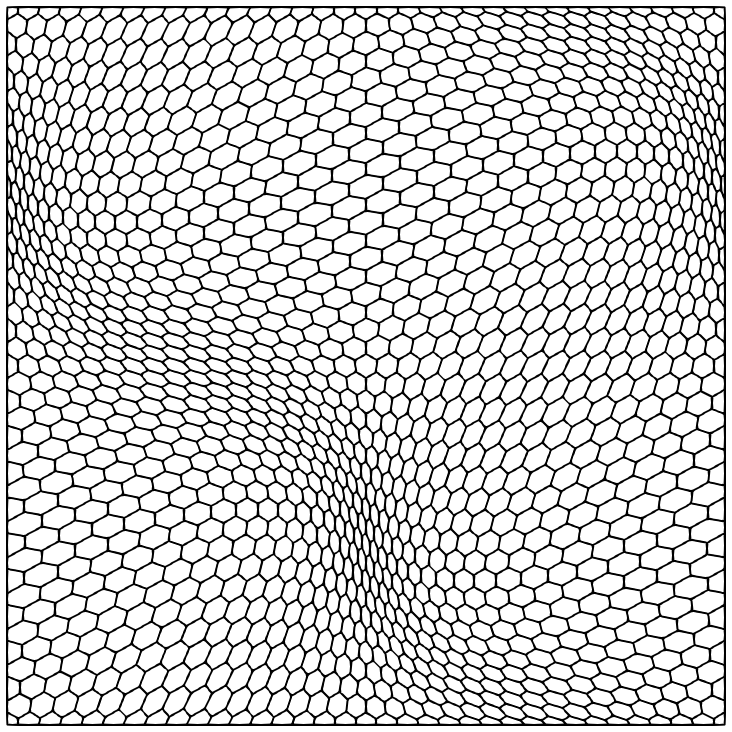}
    \caption{}
    \end{subfigure}
    \caption{Test 2: Last refinements related to the two families of meshes.}
    \label{fig:test2:mesh}
\end{figure}

We further consider two external body forces, compatible with the following two analytical solutions:
\begin{enumerate}[label=\textbf{Case \arabic*}]
    \item \label{test2:case1}\!\!\!: $\uu(\xx) = x_1(1-x_1)x_2(1-x_2) \begin{bmatrix}
        1\\
        1
    \end{bmatrix}$;
    \item \label{test2:case2}\!\!\!: $\uu(\xx) = 80 x_1(1-x_1)x_2(1-x_2) \begin{bmatrix}
        1\\
        1
    \end{bmatrix}$.
\end{enumerate}
In particular, we notice that in \ref{test2:case1} the solution is characterized by deformations of moderate magnitude, whereas much larger deformations occur in \ref{test2:case2}. For the current experiment, we consider two families of meshes. The former comprises four square meshes of $4\times4$, $8\times 8$, $16\times 16$, and $32\times 32$ square elements each. The second family is made up of $4$ regular Voronoi meshes subjected to sine distortions realized with the MATLAB library mVEM \cite{mvem}. The last refinement for each family of meshes is shown in Figure \ref{fig:test2:mesh}. Preliminary, we remark that the neural networks used to predict the NAVEM basis functions are the ones used to make the experiment reported in Section \ref{sec:nn_training}.

Due to non-linearity and possible large deformations, we here decide to test the NAVEM and VEM performance with and without the incremental force method described in Section \ref{sec:vem_discretization}, along with the Newton-Raphson scheme. Specifically, we set $N = 20$ for the incremental procedure and, for each $n=1,\dots,N-1$, we stop the Newton-Raphson scheme after $20$ non-linear iteration steps. At the incremental step $n = N$, we allow the Newton-Raphson scheme to perform all the non-linear iteration steps needed to reach the desired accuracy, i.e. we require that the relative norm of the residual is below $10^{-10}$.

We observe that, when the incremental method is applied, the VEM stabilization is chosen as the norm-based stabilization described in Section \ref{sec:vem_discretization}, where $\ww_h = \uu_h^{n-1}$. On the other hand, when we apply the plain Newton-Raphson method, without the incremental force approach, we set $\ww_h = \uu_h^0 = \bm{0}$. In particular, we observe that this last choice coincides with the ``fixed scaling'' described in \cite{DaVeiga2015}.

\begin{figure}[!ht]
    \centering
    \begin{subfigure}{0.35\linewidth}
    \includegraphics[width=1\linewidth]{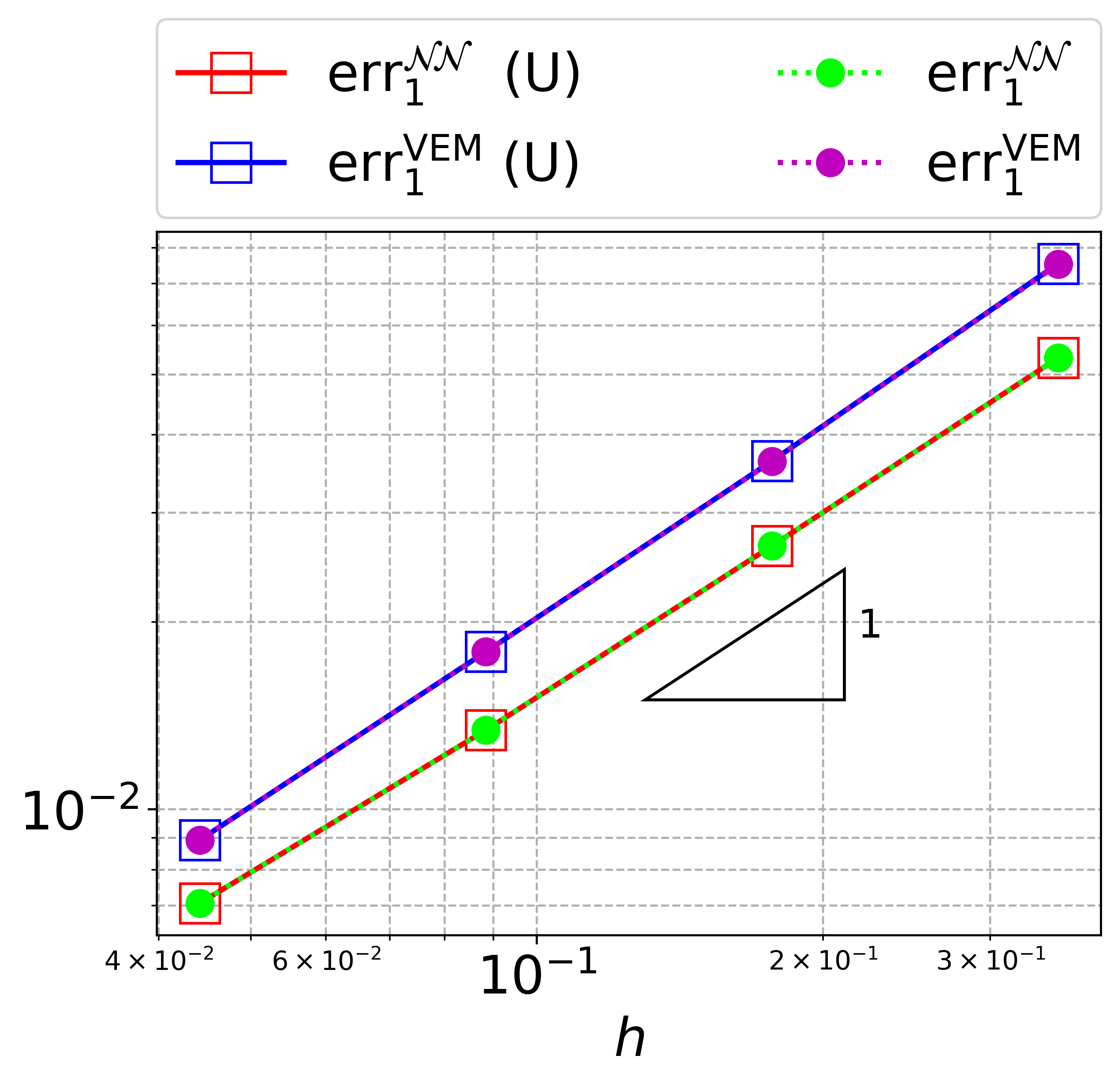}
    \caption{}
    \end{subfigure}    \hspace{30pt}  
    \begin{subfigure}{0.35\linewidth}
    \includegraphics[width=1\linewidth]{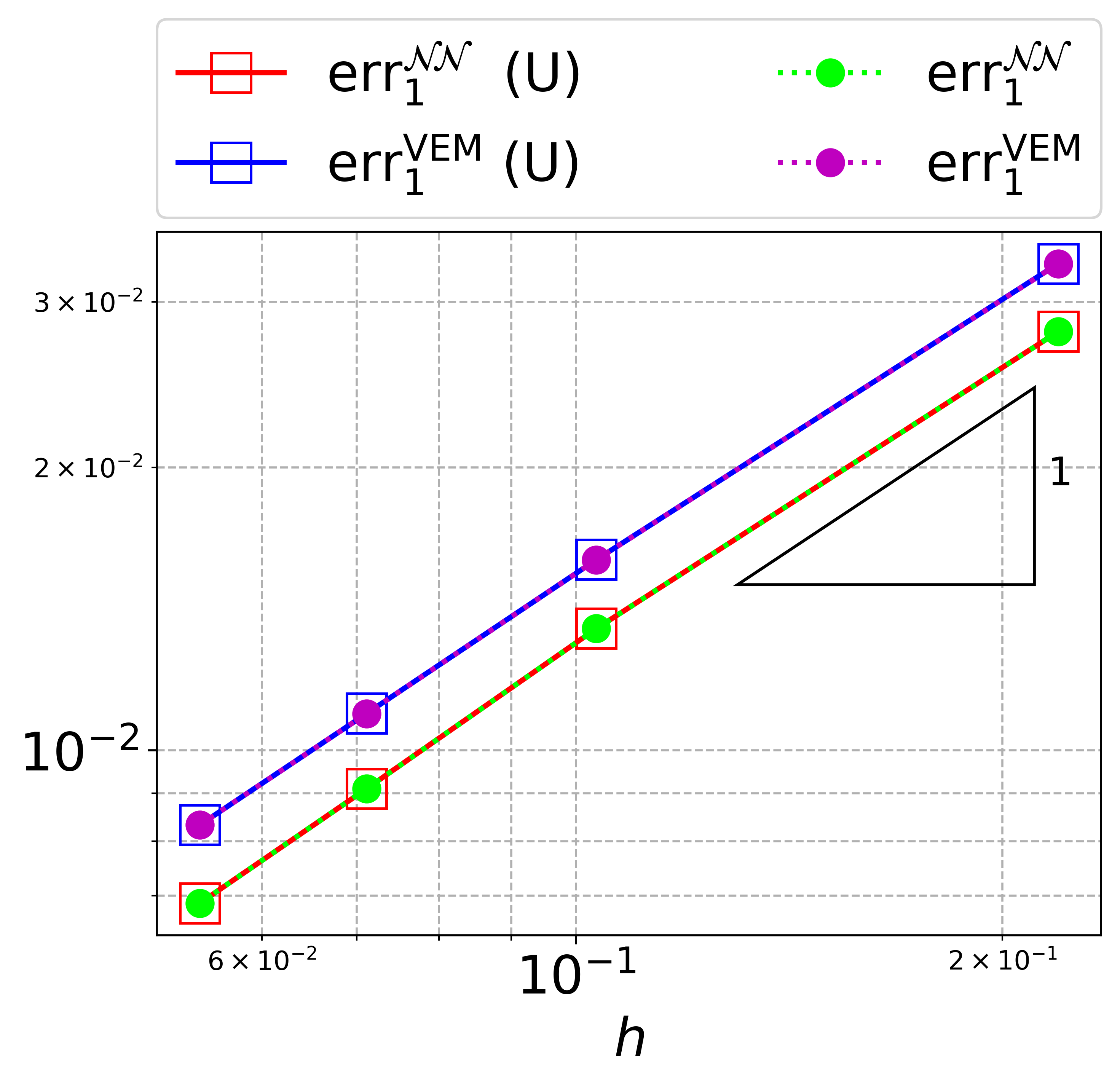}
    \caption{}
    \end{subfigure}
    \caption{Test 2, \ref{test2:case1}: Behaviour of the $H^1$-errors \eqref{eq:nn_errors}--\eqref{eq:vem_errors} as the mesh parameter $h$ decreases. ``(U)'' means that the incremental force method is adopted. Left: Cartesian grids family. Right: Distorted-Voronoi grids family.}
    \label{fig:test2:error:square}
\end{figure}

Figure \ref{fig:test2:error:square} reports the behaviour of the $H^1$-errors in \eqref{eq:nn_errors} and \eqref{eq:vem_errors}, as the mesh parameter $h$ decreases, for the test \ref{test2:case1}. Concerning this test case, we observe that the expected order of convergence is attained for both the methods with and without employing the incremental strategy. Specifically, we note that the choice of stabilization does not seem to influence the accuracy of the VEM method for both the families of meshes. Moreover, we observe that the NAVEM error curves are downward shifted with respect to the VEM error curves, as usual. Moreover, we highlight that the $L^2$-error curves follow a very similar behaviour. Thus we decide to not report this verbose information.

Figure \ref{fig:test2:error:square} reports the behaviour of $H^1$-error \eqref{eq:nn_errors}--\eqref{eq:vem_errors} as the mesh parameter $h$ decreases for the test \ref{test2:case2}, where larger displacements are recorded. Unlike the previous test case, a good choice of the stabilization term is essential to obtain good quality results when using the VEM method. Coherently with the results shown in \cite{DaVeiga2015}, we observe that adopting the incremental strategy is mandatory for VEM to avoid to level-off the stabilization term with respect to the consistency term, loosing accuracy, when Voronoi meshes are employed. On the other hand, we observe that this behaviour is less stressed in the presence of quadrilateral elements, where the stabilization term is generally less pronounced.
Concerning the NAVEM method, we observe that it does not need at all the incremental strategy, even in the presence of the larger displacement related to \ref{test2:case2}. Indeed, it behaves in the right way independently of the family of meshes and of the magnitude of the displacement for this test, while again showing a lower error constants with respect to VEM.

\begin{figure}[!ht]
    \centering
    \begin{subfigure}{0.35\linewidth}
    \includegraphics[width=1\linewidth]{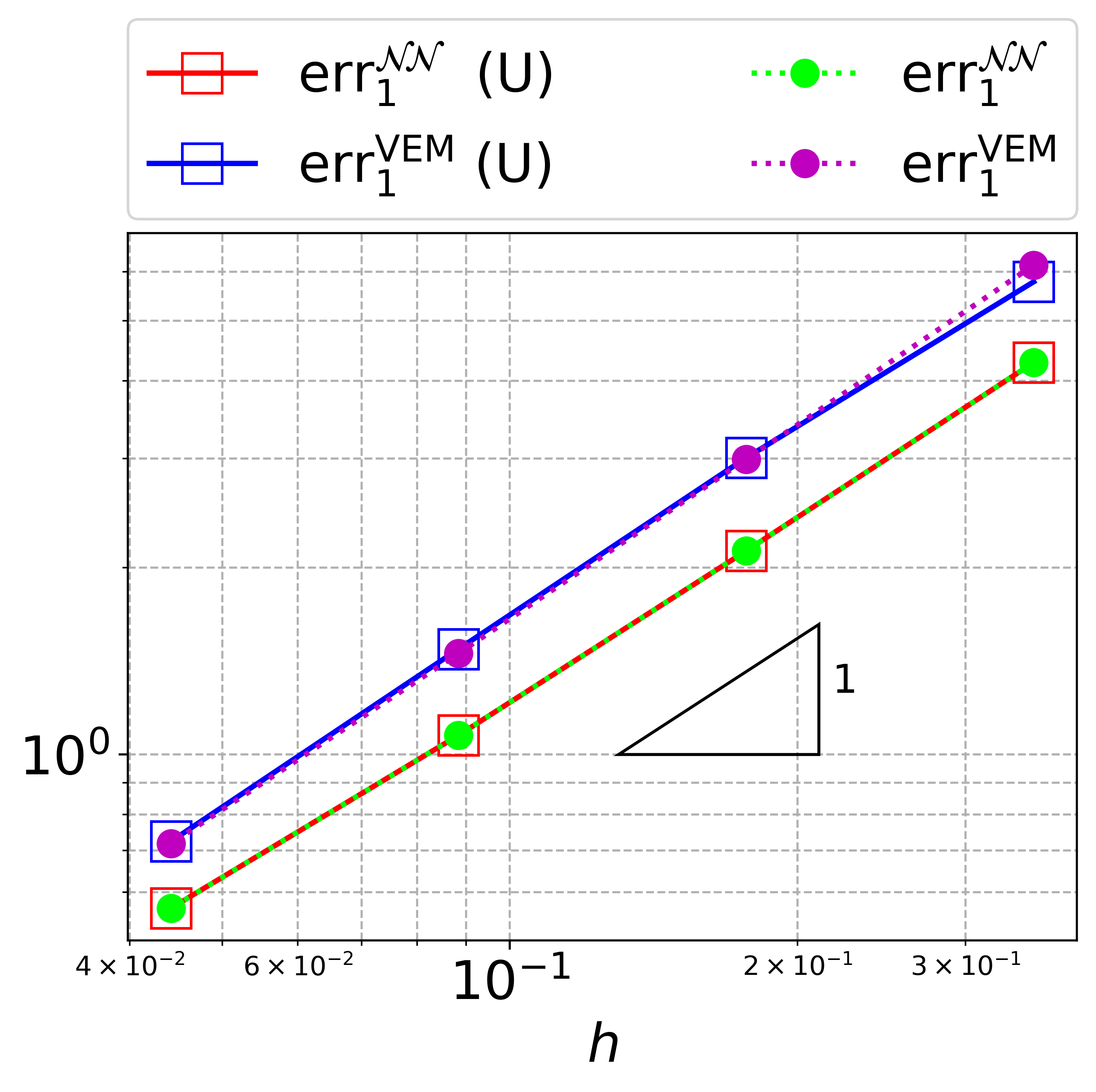}
    \caption{}
    \end{subfigure}    \hspace{30pt}  
    \begin{subfigure}{0.35\linewidth}
    \includegraphics[width=1\linewidth]{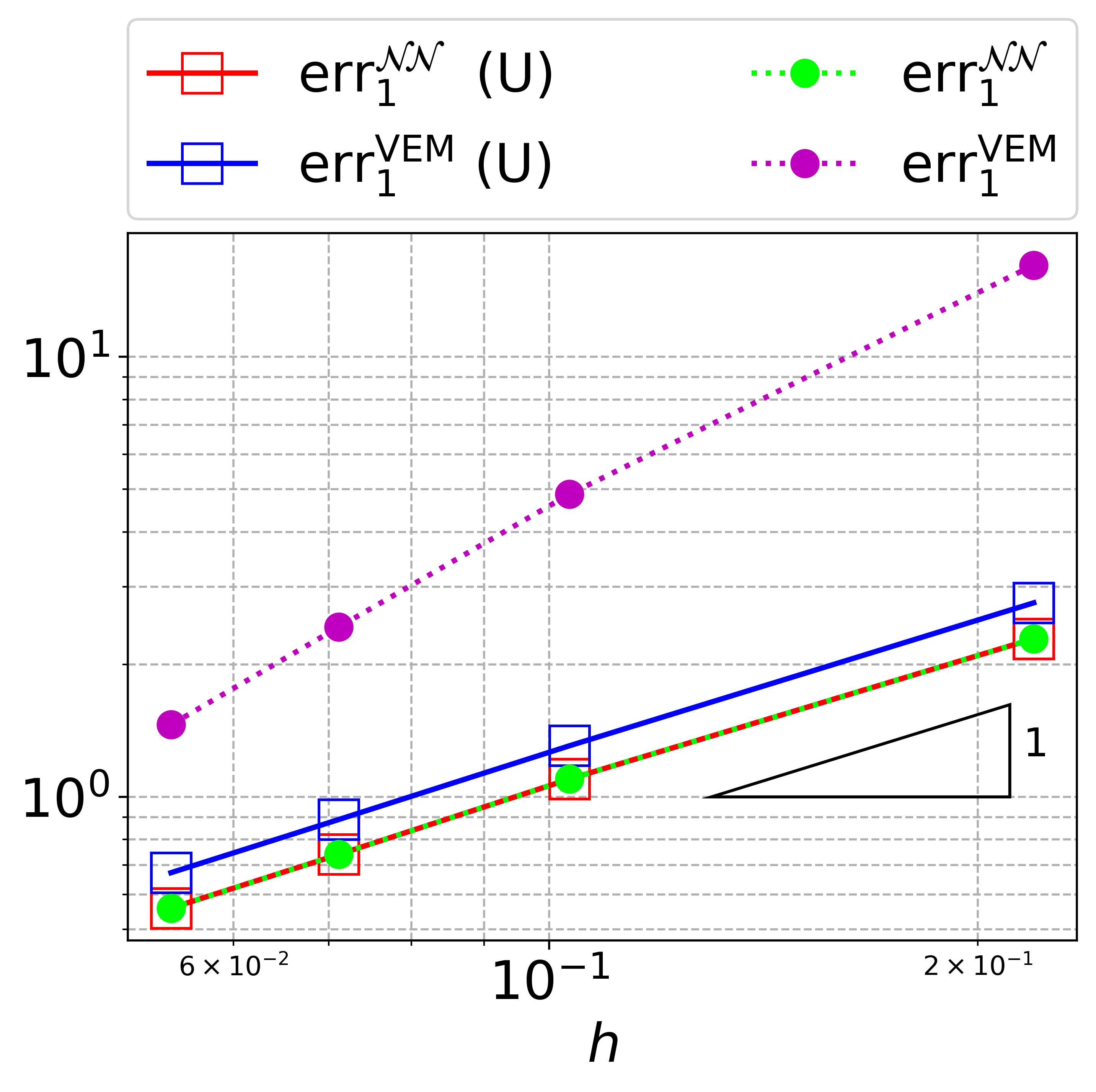}
    \caption{}
    \end{subfigure}
    \caption{Test 2, \ref{test2:case2}: Behaviour of the $H^1$-errors \eqref{eq:nn_errors}--\eqref{eq:vem_errors} as the mesh parameter $h$ decreases. ``(U)'' means that the incremental force method is adopted. Left: Cartesian grids family. Right: Distorted-Voronoi grids family.}
    \label{fig:test2:error:voronoi}
\end{figure}

\begin{figure}[!ht]
    \centering
    \includegraphics[width=0.5\linewidth]{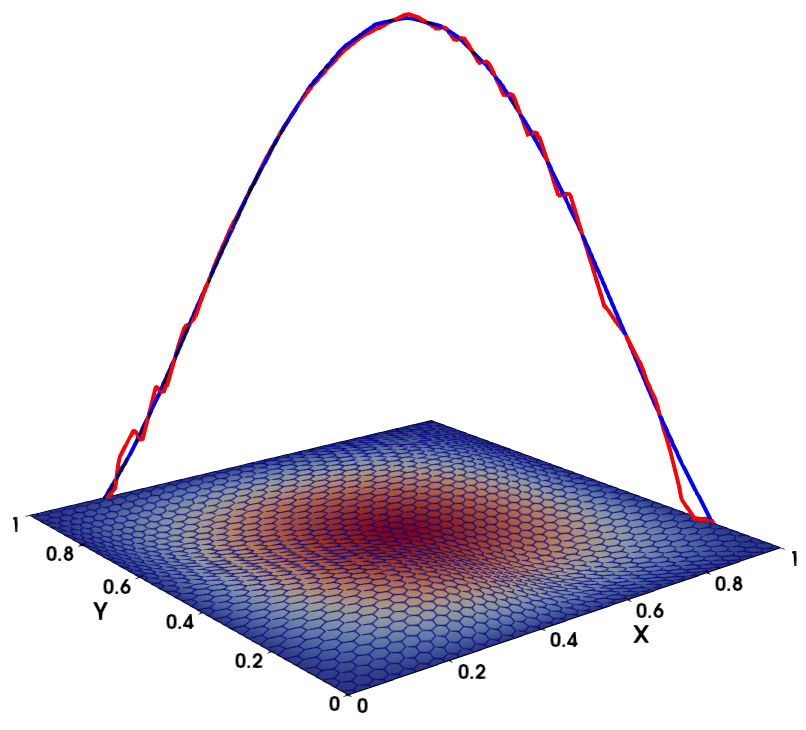}
    \caption{Test 2: Slice of the component $u_1$ related to \ref{test2:case2} and the last refinement of the Distorted-Voronoi family. Red: VEM. Blue: NAVEM. Black: exact displacement. We observe that the blue curve overlaps the black one.}
    \label{fig:test2:slice}
\end{figure}

Moreover, in Figure \ref{fig:test2:slice}, we observe that NAVEM does not suffer of the same unstable spurious oscillating behaviour, remarked also in \cite{DaVeiga2015}, that characterizes the VEM solution when a fixed scaling is employed for the \ref{test2:case2}. In this figure, we draw a slice of the first component of the VEM solution, of the NAVEM solution, and of the exact solution. It is evident that the VEM solution is characterized by large oscillations that ruin its accuracy, whereas the other two curves are almost exactly overlapped and are indistinguishable in the plot.

Finally, we observe that the robustness that characterizes NAVEM translates to a higher computational efficiency with respect to VEM. Indeed, when the incremental strategy is not employed, the Newton-Raphson method requires $125$ non-linear steps to reach the desired accuracy for the last refinement of Distorted-Voronoi family and \ref{test2:case2}. Whereas, since VEM requires the incremental strategy to perform a correct choice for the stabilization, the whole non-linear incremental strategy needs about $500$ steps to converge.

\subsection{Test 3: Neo-Hookean Hyperlasticity}

In this last test case, we consider Problem \eqref{eq:model_problem} with $\ssigma$ corresponding to the non-symmetric first Piola-Kirchhoff stress tensor $\Pmatrix$. More precisely, we focus on an hyperelastic material of Neo-Hookean type by setting
\begin{equation}
    \Pmatrix(\Fmatrix) = \mu(\Fmatrix - \Fmatrix^{-T}) + \lambda J \Theta(J) \Theta'(J) \Fmatrix^{-T},
    \label{eq:piola_neo_hook}
\end{equation}
where $\Fmatrix$ is the deformation gradient $\Fmatrix = \Imatrix + \nabla \uu$, $J \in \R_{>0}$ is the determinant of $\Fmatrix$, and $\Theta: \R_{>0} \to \R$ is a smooth function such that $\Theta(J) = 0\ \Leftrightarrow \ J = 1$ and $\Theta'(1) \neq 0$ \cite{AbbasErn2018}. In equations \eqref{eq:piola_neo_hook}, $\lambda$ and $\mu$ are given positive constants.

We observe that, to carry out the discretization related to a Neo-Hookean material, a discrete version of the determinant $J$ must be provided. In the Virtual Element framework, two possible discretization are considered. The simplest choice coincides with taking on each element $E \in \Th$, for each $k \geq 0$ and $n \geq 1$:
\begin{equation*}
    J = \det (\Imatrix + \nabla \uu_h^{VEM, n, k}) \approx \det( \Imatrix + \proj{0,E}{0} \nabla \uu_h^{VEM, n, k}).
\end{equation*}
Again, since the virtual functions are not known in a closed form in the interior of the element, a projector operator should be introduced to access the point-wise evaluations of such functions. With this kind of approximation, the determinant $J$ is constant on each element $E \in \Th$. A more stable version is proposed in \cite{Chi2016}, where $J$ is approximated with its mean value over $E \in \Th$. More precisely, the authors in \cite{Chi2016} approximate $J$ as follows:
\begin{equation}
    J \approx \frac{1}{\vert E \vert} \int_{E} J = \int_{\widetilde{E}} ~d \widetilde{\xx} = \frac{\vert \widetilde{E} \vert}{\vert E \vert},
    \label{eq:stable_det_vem}
\end{equation}
where $\vert E \vert$ denotes the measure of the element $E \in \Th$, whereas $\widetilde{\cdot}$ defines a quantity $\cdot$ in the current deformed configuration. 
In the NAVEM framework, we evaluate the determinant simply as
\begin{equation}
     J = \det (\Imatrix + \nabla \uu_h^{\NN, n, k}),
\end{equation}
since we are able to evaluate NAVEM functions and their gradients everywhere. 

We test the method considering 
\begin{equation*}
    \Theta(J) = \log J,\quad \Theta'(J) = \frac{1}{J},\quad \lambda = 5.1,\quad \mu = 1.0,
\end{equation*}
whose values of Lamé coefficients correspond to a Poisson ratio $\approx 0.45$. We observe that in this case, the elastic modulus $\mathbb{A}$ in \eqref{eq:elastic_modulus} in this specific case depends only on the deformation gradient and it is given by \cite{AbbasErn2018}:
\begin{align*}
    \mathbb{A}(\Fmatrix) = \frac{\partial \Pmatrix(\Fmatrix)}{\partial \Fmatrix} &= \mu (\Imatrix\ \overline{\otimes}\ \Imatrix + \Fmatrix^{-T} \underline{\otimes} \Fmatrix^{-1}) - \lambda J \Theta(J) \Theta'(J) \Fmatrix^{-T} \underline{\otimes} \Fmatrix^{-1}\\
    &\quad + \lambda J \left[\Theta(J)\Theta'(J) + J (\Theta'(J))^2 + J \Theta(J) \Theta''(J)\right] \Fmatrix^{-T} \otimes \Fmatrix^{-T},
\end{align*}
where $\otimes,\ \underline{\otimes}$ and $\overline{\otimes}$ are defined such that $\{\Amatrix \otimes \Bmatrix\}_{ijkl} = \Amatrix_{ij} \Bmatrix_{kl}$, $\{\Amatrix \underline{\otimes} \Bmatrix\}_{ijkl} = \Amatrix_{il} \Bmatrix_{jk}$  and $\{\Amatrix \overline{\otimes} \Bmatrix\}_{ijkl} = \Amatrix_{ik} \Bmatrix_{jl}$, for all $i,j,k,l = 1,\dots, 2$ and for each pair of second-order tensors $\Amatrix$ and $\Bmatrix$.

The initial configuration $\Omega$ is set as the unit square $[0,1]^2$ and the displacement is clamped at $\Gamma_D = \{0\} \times [0,1]$, while the remaining part of the boundary is free. The body load is given by
\begin{equation}
    \ff(x_1, x_2) = [100 x_2^3, 0] \quad \forall (x_1,x_2) \in \Omega.
\end{equation}

\begin{figure}[!ht]
    \centering
    \begin{subfigure}{0.3\linewidth}
    \includegraphics[width=1\linewidth]{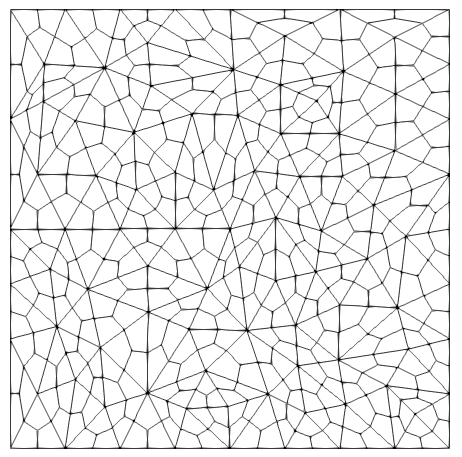}
    \caption{}
    \label{fig:test3:mesh}
    \end{subfigure}   
    \begin{subfigure}{0.33\linewidth}
    \includegraphics[width=1\linewidth]{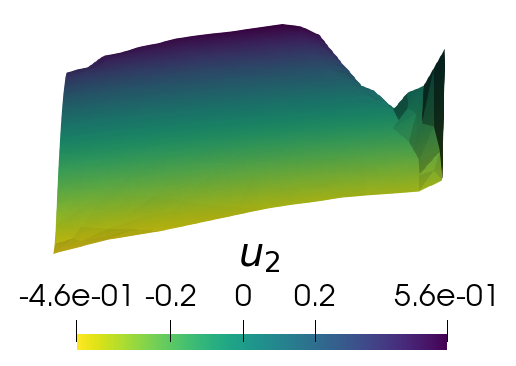}
    \caption{}
    \label{fig:test3:navem_sol}
    \end{subfigure} 
    \begin{subfigure}{0.33\linewidth}
    \includegraphics[width=1\linewidth]{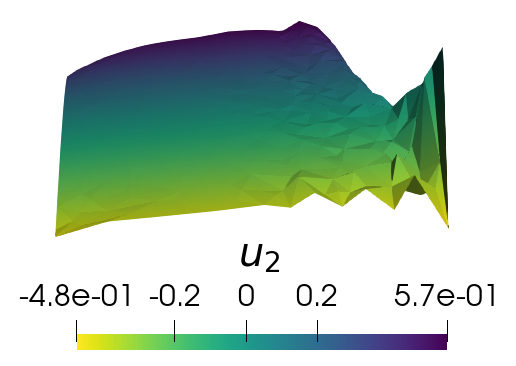}
    \caption{}
    \label{fig:test3:vem_sol}
    \end{subfigure}
    \caption{Test 3: Left: Mesh used in the simulation. Center: Displacement in the $x_2$-direction obtained with the NAVEM method. Right: Displacement in the $x_2$-direction obtained with the VEM method.}
    \label{fig:test3:solution}
\end{figure}

\begin{figure}[!ht]
    \centering  
    \begin{subfigure}{0.3\linewidth}
    \includegraphics[width=1\linewidth]{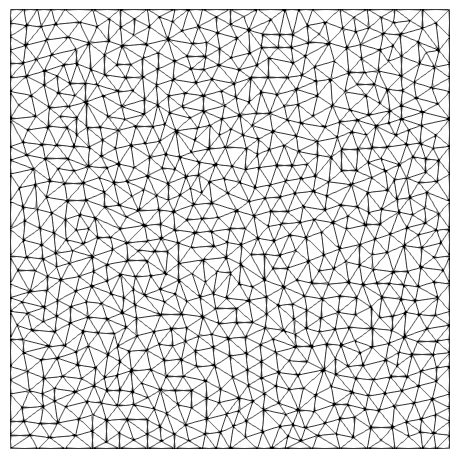}
    \caption{}
    \label{fig:test3:ref_mesh}
    \end{subfigure} \hspace{30pt}
    \begin{subfigure}{0.33\linewidth}
    \includegraphics[width=1\linewidth]{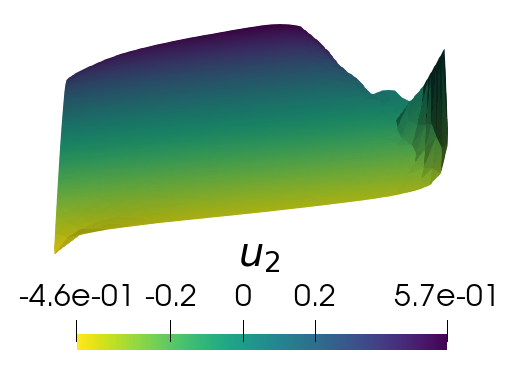}
    \caption{}
    \label{fig:test3:fem_sol}
    \end{subfigure}
    \caption{Test 3: Left: Mesh used to compute the FEM reference solution. Right: Displacement in the $x_2$-direction obtained with the FEM method.}
    \label{fig:test3:fem_solution}
\end{figure}

We simulate this test on the mesh shown in Figure \ref{fig:test3:mesh}, generated as in \cite{Lamperti2022}, performing $N = 50$ iterations of the incremental force method. At each step $n = 1,\dots,N$, we iterate the Newton-Raphson method until the residual norm is below $1.0e-09$. Figures \ref{fig:test3:navem_sol} and \ref{fig:test3:vem_sol} show the $x_2$-component of the discrete displacement obtained with NAVEM and VEM, respectively, whereas Figure \ref{fig:test3:fem_sol} shows the result of the Finite Element discretization of the problem on the finer triangular mesh \ref{fig:test3:ref_mesh}. 
We first observe that, even though we employ the more stable approximation of determinant, defined in \eqref{eq:stable_det_vem}, to produce virtual element results, we do not appreciate a significant difference between the two approximation procedures for the determinant in this test case.
Moreover, we highlight that all the aforementioned methods take about $5$ Newton steps for each incremental step $n$ to reach the desired accuracy. Despite this, we can observe some spurious oscillations in the $x_2$-component of the virtual element displacement that are not present in the NAVEM or FEM solution. 


To conclude, we observe that the NAVEM has a simpler formulation than the VEM one and, in the provided numerical tests, it is more stable in the presence of large deformations. On the other hand, despite having a formulation similar to the FEM one, it ensures greater flexibility by allowing the use of more general meshes.

\section{Conclusions}\label{sec:conclusion} 
In this paper, we present the Neural Approximated Virtual Element Method to solve linear and non-linear elasticity problems. The NAVEM is a polygonal method that approximates via neural networks the function spaces used in the Virtual Element Method. This neural approximation provides closed-form local basis functions, eliminating the need for additional stabilization and projection operators. 

The absence of a stabilization term is particularly advantageous in nonlinear problems, where designing a suitable and computable stabilization operator that preserves accuracy can be challenging. Elasticity problems represent a particular domain in which numerous different non-linearities may be present, depending on the simulated materials. For this class of problems, we discuss the NAVEM formulation, emphasizing the simpler form with respect to the standard VEM formulation. Numerical results are in agreement with what has been observed for elliptic problems in \cite{PintoreTeora2024, PintoreTeora2025} and highlight the benefits of bypassing stabilization operators.

Future perspectives include extending NAVEM to more physically relevant problems, such as handling internal constraints like incompressibility. Additionally, crack propagation and adaptivity present promising applications for this method. Indeed, as observed in \cite{PintoreTeora2025}, the method can be easily optimized for mesh comprising triangles with hanging nodes. This feature opens the possibility of integrating NAVEM with standard finite element methods by refining only elements where non-linearities are most intense: introducing hanging nodes instead of globally remeshing neighboring elements.

\section*{Acknowledgements}
The author S.B. kindly acknowledges partial financial support provided by PRIN project ``Advanced polyhedral discretisations of heterogeneous PDEs for multiphysics problems'' (No. 20204LN5N5\_003), by PNRR M4C2 project of CN00000013 National Centre for HPC, Big Data and Quantum Computing (HPC) (CUP: E13C22000990001) and the funding by the European Union through project Next Generation EU, M4C2, PRIN 2022 PNRR project P2022BH5CB\_001 ``Polyhedral Galerkin methods for engineering applications to improve disaster risk forecast and management: stabilization-free operator-preserving methods and optimal stabilization methods.''. The author M.P. kindly acknowledges financial support provided by PEPR/IA (\url{https://www.pepr-ia.fr/}). The author G.T. kindly acknowledges the financial support provided by INdAM-GNCS Project ``Metodi numerici efficienti per problemi accoppiati in sistemi complessi'' (CUP: E53C24001950001). 


\bibliographystyle{elsarticle-num} 
\bibliography{biblio.bib}



%
%
%
\end{document}